\documentstyle[art10]{article}
\title{Classification of \\ Plane Congruences in 
 {\bf P}$^4_{\bf C}$ (II)}
\author{Manuel Pedreira 
\and Luis-Eduardo Sol\'{a}-Conde\thanks{Supported by an F.P.U.
fellowship of Spanish Government}}
\date{}
\newtheorem{teorema}{Theorem}[section]

\newtheorem{corolario}{Corollary}[teorema]

\newtheorem{nota}{Remark}[section]

\def\rank{\mathop{\rm rank}}

\begin{document}
\maketitle

{\footnotesize{\bf Authors' address:} Departamento de Algebra, Universidad de Santiago
de Compostela. $15706$ Santiago de Compostela. Galicia. Spain. e-mail: {\tt
pedreira@zmat.usc.es}\\
{\bf Abstract:} We study second order focal loci of nondegenerate plane congruences in
${\bf P}^4$ with  degenerate focal conic. We show the projective
generation of such congruences when the second order focal locus fills a component
of the focal conic, improving a result of Corrado Segre.\\
{\bf Mathematics Subject Classifications (1991):} Primary, 14N05; secondary,
51N35.\\
{\bf Key Words:} Focal points, congruences.}

\vspace{0.1cm}

{\bf Introduction:} Throughout this paper, the base field for algebraic
varieties will be ${\bf C}$. Let ${\bf P}^4$ be the $4$--dimensional complex
projective space and $G(2,4)$ the Grassmannian of planes in ${\bf P}^4$.
Given a quasi-projective irreducible algebraic surface $\Sigma\subset G(2,4)$, 
let ${\bf P}^2(\sigma)$ be the plane in ${\bf P}^4$ corresponding to each 
$\sigma\in\Sigma$ and let ${\cal I}_{\Sigma}:=\{(P,\sigma)\in{\bf
P}^4\times\Sigma:\,P\in{\bf P}^2(\sigma)\}$ denote the Incidence Var\-iety.
Consider the two projections $p:{\cal I}_{\Sigma}\longrightarrow{\bf P}^4$ and
$q:{\cal I}_{\Sigma}\longrightarrow\Sigma$; we call $(\Sigma,p,q)$ a plane
congruence with projective realization $V(\Sigma):=\overline{p({\cal
I}_{\Sigma})}$. If $V(\Sigma)={\bf P}^4$, the congruence is said to be
nondegenerate. From now on, a plane congruence will be nondegenerate. Recall that the
first order focal locus is defined by
$F_1(\Sigma):=\{(Q,\sigma)\in{\cal I}_{\Sigma}:\,(dp)_{(Q,\sigma)}\mbox{ is not
injective }\}$, which is a proper closed subscheme of ${\cal I}_{\Sigma}$. There is a
nonempty open set
$U\subset\Sigma$ consisting of smooth points such that $F_1(\Sigma)\cap{\bf
P}^2(\sigma)=C(\sigma)$ is a conic. In [3] we have classified the plane
congruences in ${\bf P}^4$ such that $C(\sigma)$ is degenerate. We  will refer to an
$\alpha$,$\beta$, $\gamma$ or $\delta$--Congruence in according to the existence of
one, two (different or coincident), or an infinity of developable families passing
through the  general plane ${\bf P}^2(\sigma)$ of the congruence. It's well known, see
[2]
$1.5$, that the second order foci of a congruence are either at most $5$ points or
they fill a component of
$C(\sigma)$. In relation with our classification, we show in $\S 1$ that there is only
one  second order focal point for the general $\gamma$--congruence and there are 
just three second order focal points for the general $\alpha$ or
$\beta$--congruence. If the second order foci fill a component of $C(\sigma)$, we
obtain either a $\delta$--congruence or degenerate cases of $\alpha$, $\beta$, $\gamma$
that we  describe in $\S 2$. Our analysis includes the proof of a footnote assertion
due to Corrado Segre in [4] which appeared with a partial statement in [2] $1.6$.

\smallskip

The results on this paper belong to the Ph.D. thesis of the second author with the
advice of the first one. We thank Ciro Ciliberto by reading a first version of
this paper and his criticism.

\addtocontents{toc}{\protect\vspace{3ex}}
\section{Second order focal locus for nondegener\-ate plane congruences with
de\-gener\-ate fo\-cal conic.}
\addtocontents{toc}{\protect\vspace{1ex}}

Given a congruence $(\Sigma, p, q)$, we consider the family of focal conics ${\bar
q}:F_1(U)\longrightarrow U$ and the corresponding projective realization ${\bar
p}:F_1(U)\longrightarrow{\bf P}^4$. The second order focal locus of the congruence
$(\Sigma,p,q)$ is defined to be the first order focal locus of $(U,{\bar p},{\bar q})$. We
refer to [1] for details about higher order focal locus. To compute the first order focal
locus corresponding to $(U,{\bar p},{\bar q})$, we have to take a desingularization of 
$F_1(U):=\{(Q,\sigma)\in{\cal I}_U:\,\rank(dp)_{(Q,\sigma)}\leq 3\}$. There are
three possibilities:
\begin{enumerate}
\item If $C(\sigma)$ is nondegenerate for $\sigma\in U$, $F_1(U)$ is smooth. To
compute the focal locus of $(U,{\bar p},{\bar q})$ on the
conic $C(\sigma)$, we consider the characteristic map at $\sigma\in U$,
$\lambda:T_{U,\sigma}\longrightarrow H^0(C(\sigma),{\cal N}_{C(\sigma),P^4})$. If
$T_{U,\sigma}=\langle v_1,v_2\rangle$, since ${\cal N}_{C(\sigma),P^4}=({\cal
O}_{C(\sigma)}(1))^2\oplus{\cal O}_{C(\sigma)}(2)$, then 
 $\lambda(v_1)=(f_{11},f_{12},F_{13})$, $\lambda(v_2)=(f_{21},f_{22},F_{23})$.
When $(\lambda:\mu)\in{\bf P}^1$, the solution of the system $\lambda f_{11}+\mu
f_{21}=0$; $\lambda f_{12}+\mu f_{22}=0$ is the first order focal point
$P(\lambda:\mu)\in C(\sigma)$ corresponding to the direction $(\lambda:\mu)$.
So $C(\sigma)$ is the image of the Veronese's map $(\lambda:\mu)\in{\bf
P}^1\longrightarrow P(\lambda:\mu)\in {\bf P}^2(\sigma)$. The point
$P(\lambda:\mu)$ will be a focal point of the family of conics if $\lambda
F_{13}(P(\lambda:\mu))+\mu F_{23}(P(\lambda:\mu))=0$. Then the second order foci are
either five points on $C(\sigma)$ or the whole conic.
\item If $C(\sigma)=(r(\sigma))^2$ is a double line, $F_1(U)$ is smooth
consisting of a $2$--dimensional family of lines. The second order focal locus
of the congruence $(\Sigma,p,q)$ is defined by $F_2(U):=\{(Q,\sigma)\in
F_1(U):\,\rank(d{\bar p})_{(Q,\sigma)}\leq 2\}$.
\item If $C(\sigma)=r(\sigma)\vee r'(\sigma)$ consists of two different lines,
we consider the subvarieties $R(U):=\{(Q,\sigma)\in{\cal I}_U:\,Q\in
r(\sigma)\}$ and $R'(U):=\{(Q,\sigma)\in{\cal I}_U:\,Q\in r'(\sigma)\}$. We get a
desingularization by taking $u:R(U)\sqcup R'(U)\longrightarrow F_1(U)$ . In this case, a
second order focal point for the congruence
$(\Sigma,p,q)$ is a first order focal point for any of the families of lines
${\bar q}_1:R(U)\longrightarrow U$ or ${\bar q}_2:R'(U)\longrightarrow U$.
\end{enumerate}

\begin{nota}{\rm Observe, for example, that given the family ${\bar
q}_1:R(U)\longrightarrow U$ and the characteristic map at $\sigma\in U$,
$\lambda:T_{U,\sigma}\longrightarrow H^0(r(\sigma),{\cal N}_{r(\sigma),P^4})$, since ${\cal
N}_{r(\sigma),P^4}\cong {\cal N}_{r(\sigma),P^2(\sigma)}\oplus({\cal
N}_{P^2(\sigma),P^4})|_{r(\sigma)}$, $\lambda(w)=(\psi_r(w),\chi(w))$, with $\chi$ the
restriction of the characteristic map at $\sigma$ corresponding to the congruence
$(\Sigma,p,q)$. So $Q\in r(\sigma)$ is a second order focal point if there is $w\in
T_{\Sigma,\sigma}$ such that
$\chi(w)(Q)=\psi_r(w)(Q)=0$. That is, there is $w\in\ker(dp)_{(Q,\sigma)}$ such that
$\psi_r(w)(Q)=0$.} 
\end{nota}

We compute the second order foci for the plane congruences classified on [3]:

\smallskip

$\delta$--{\bf Congruences:} In this case $(\Sigma,p,q)$ consists of the planes
containing a line $r\subset{\bf P}^4$ by [3], Theorem $2.1.1$. The first order
focal locus on every plane is the double line $C(\sigma)=r^2$ being $r={\bar
p}(F_1(U))$. Every point of $r$ is a fundamental point and so, it's a second
order focal point by [2], Lemma $1.7$.\\
$\beta$--{\bf Congruences:} The focal conic consists of two different lines
$r(\sigma)$ and $r'(\sigma)$, which are the focal loci corresponding to two
directions $v$ and $v'$ in $T_{\Sigma,\sigma}$, and the focal loci for the remaining
directions consists of the point $P(\sigma)=r(\sigma)\cap r'(\sigma)$. We use 
Remark $1.1$ to compute the second order focal points of the congruence on one of the
focal lines, for example $r(\sigma)$. The second order focal locus for the developable
direction $v\in T_{\Sigma,\sigma}$ consists of the points $Q\in r(\sigma)$ such
that $\psi_r(v)(Q)=0$; that is, it is either one point or the whole
line, in correspondence with the developable system passing through ${\bf P}^2(\sigma)$
for the direction $v$. For the remaining directions $\lambda v+v'$, the only first order
focal point is $P(\sigma)$, and the relation $\psi_r(\lambda v+v')(P(\sigma))=0$ shows
that $P(\sigma)$ is a second order focal point for one of these directions.\\
$\gamma$--{\bf Congruences:} The focal conic consists of a double line
$C(\sigma)=(r(\sigma))^2$ being $r(\sigma)$ the focal locus for one (double)
direction $v$. The focal loci for the remaining directions consist of one point $P(\sigma)\in
r(\sigma)$. As in case $\beta$, $P(\sigma)$ will be a second order focal point for some
direction. For the developable direction $v$, the second order focal locus 
is given by the solutions of a linear form  $\psi_r(v)=0$; that is, either one
point or the whole line. If it consists of one point, this point must be $P(\sigma)$. In
fact, recall from [3] the classification of $\gamma$--congruences given in 2.2 and use
proposition 1.1.1: If the developable families contained in $\Sigma$ are not families of
planes containing a line, their second order focal points must lie on the variety
$\Sigma'$ that is the projective realization of the points $P(\sigma)$, $\sigma\in U$.\\
$\alpha$--{\bf Congruences:} Suppose $C(\sigma)=r(\sigma)\vee r'(\sigma)$ being $r(\sigma)$ the
focal line for one direction $v\in T_{\Sigma,\sigma}$. For the other directions
$v'+\lambda v$ there is only one focal point on $r'(\sigma)\setminus r(\sigma)$. The
second order focal points on $r(\sigma)$ correspond to the solutions of the
linear equation $\psi_r(v)=0$, and so to either only one point or the whole line,
corresponding to the developable system passing through ${\bf P}^2(\sigma)$. The second
order focal locus on $r'(\sigma)$ corresponds to the points $Q\in r'(\sigma)$ such
that $\chi(\lambda v+\mu v')(Q)=\psi_{r'}(\lambda v+\mu v')(Q)=0$ for some
$(\lambda:\mu)$. Equivalently, 
$\det\left(^{\chi(v)(Q)\quad\chi(v')(Q)}_{\psi_{r'}(v)(Q)\,\,\psi_{r'}(v')(Q)}\right)=0$,
then it consists of either two points on $r'(\sigma)$ or the whole line.

\addtocontents{toc}{\protect\vspace{3ex}}
\section{Nondegenerate plane congruences with second order focal locus filling
a component of the focal conic.}
\addtocontents{toc}{\protect\vspace{1ex}}

Suppose $(\Sigma,p,q)$ is a nondegenerate plane congruence in ${\bf P}^4$ with 
general focal conic $C(\sigma)=r(\sigma)\vee r'(\sigma)$ such that $r(\sigma)$ is the focal
locus on ${\bf P}^2(\sigma)$ for one direction. We have seen in the above section that 
$r(\sigma)$ consists of second order focal points iff the corresponding
developable system is a family of planes containing $r(\sigma)$. The projective
realization of this family of lines, $p(R(U))$, will be either a line in case
$\delta$, or a ruled surface. In accordance to the classification given in [3] we
get the following.

\begin{teorema} Let $(\Sigma,p,q)$ be a nondegenerate plane congruence in ${\bf
P}^4$ with degenerate general focal conic such that the second order focal locus
fills some of the focal lines corresponding to developable directions. Then, either
$(\Sigma,p,q)$ is a $\delta$--congruence or the projective realization of such
lines is a ruled surface $S$ such that there is a $1$--dimensional family of
planes in the congruence passing through each one of the generators of $S$. We have
the following cases:
\begin{enumerate}
\item $S$ is nondevelopable and $(\Sigma,p,q)$ consists of the tangent planes to
$S$ ($\gamma_1$--congruence).
\item $S$ is nondevelopable and $(\Sigma,p,q)$ is not the family of tangent planes to
$S$ ($\alpha_3$--congruence).
\item $S$ is developable and the $1$--dimensional families containing each one
of the generators of $S$ are linear systems containing the tangent plane to $S$ along the
generator. $S$ may be either a cone ($\gamma_3$--congruence) or a tangent developable
($\gamma_2$--congruence).
\item $S$ is developable and the $1$--dimensional families containing each
generator are not as in  case $3$ above. $S$ may be either a cone
($\beta_3$--congruence) or a tangent developable ($\beta_2$--congruence). 
\end{enumerate}
\end{teorema}

{\bf Proof:} From [3], Theorems $2.3.3$ and $2.2.2$, we obtain assertions $1$,
$2$, and $3$, $4$, respectively. We analyze the plane congruences consisting of
$1$--dimensional families of planes passing through each one of the generators of a
cone $S$ with vertex $\Sigma'$. They are $\beta_3$/$\gamma_3$--congruences. We
study the existence of another developable system passing through each plane. Suppose
$\pi\subset{\bf P}^4$ is a general plane; there is an open set $U\subset\Sigma$
such that the map $\psi:U\longrightarrow\pi$, $\psi(\sigma)={\bf P}^2(\sigma)\cap\pi$ is
well defined and dominant. Let $H\subset{\bf P}^4$ be a general hyperplane not passing by
$\Sigma'$ and consider the map $\rho:U\longrightarrow C:=S\cap H$,
$\rho(\sigma)=r(\sigma)\cap H$. For every $\sigma\in U$, we have ${\bf
P}^2(\sigma)=\langle\Sigma',\psi(\sigma),\rho(\sigma)\rangle$, and we can take a
parametrization $\sigma(u,v)$ of $U$ around $\sigma$ and a reference system on $\pi$
such that $\sigma(0,0)=\sigma$, $\psi(\sigma(u,v))=(u:v:1)$. For every family of planes
passing through the generator $\langle\Sigma',\rho(\sigma)\rangle$ of $S$, the section
by the plane $\pi$ is a curve $C_{\rho(\sigma)}$ satisfying the differential equation
$dv-kdu=0$ around $\sigma$ for some function $k(u,v)$. To conclude, we argue in the
same way that the proof of $2.2.2$ in [3]. \rule{2mm}{2mm}

\begin{teorema} Let $(\Sigma,p,q)$ be an $\alpha$--congruence such that the second order focal
locus on every plane fills the focal line not corresponding to a developable direction. The
projective realization of these lines is a plane $\pi$ and every developable
system must be contained in a hyperplane passing through $\pi$.
\end{teorema}

{\bf Proof:} Suppose $r'(\sigma)\subset C(\sigma)$ is the focal line not corresponding
to a developable direction and consisting of second order focal points. Let ${\bar
q}_2:R'(U)\longrightarrow U$ be the family of these lines. Since every point is a focal
point for this family, $\dim p(R'(U))<\dim R'(U)$. Moreover the morphism
$\psi:U\longrightarrow G(1,4)$ defined by $\psi(\sigma):=[r'(\sigma)]$ has finite fibers
because otherwise $r'(\sigma)$ would be the focal line for the developable system of
planes containing it. Then $\dim\overline{\psi(U)}=2$ and the projective realization
$\overline{p(R'(U))}$ is a surface containing a $2$--dimensional family of lines; and
thus a plane $\pi$. 

\smallskip

For the other assertion, if the developable systems consist of planes containing the line
$r(\sigma)$, every developable system is contained in the hyperplane $\langle\pi,
r(\sigma)\rangle$. So we have to prove that, if $S\subset{\bf P}^4$ is a developable
ruled surface such that their tangent planes meet a plane $\pi$ on lines, then $S$ is
degenerate and contained in a hyperplane passing by $\pi$. If $S$ is the tangent
developable to a curve $C$, every one of the tangent lines meet $\pi$, so $C$ and $S$
lie in a hyperplane
$H$. If $H$ does not contain $\pi$, every tangent line to $S$ meet the line $\pi\cap H$ and
$C$ is contained in a plane $\pi'$, so $S\subset\langle\pi,\pi'\rangle$. Now, suppose $S$
is a cone and take a directrix $C=H_0\cap S$ cutting with a hyperplane $H_0$ not passing
through the vertex $O$. The tangent lines to $C$ meet the line $\pi\cap H_0$, so $C$
lies in a plane
$\pi'$ and $S\subset\langle\pi',O\rangle=:H$. The same argument with other hyperplane
$H_1$ such that $\pi\cap H_0\neq\pi\cap H_1$ yields that $H$ contains $\pi\cap H_0$ and
$\pi\cap H_1$, and then $H\supset\pi$. \rule{2mm}{2mm}

\begin{corolario} With the same hypothesis as in the above theorem, an
$\alpha_1$--congru\-ence
$(\Sigma,p,q)$ consists of the osculating planes to a $1$--dimensional linear system of
hyperplane sections of a nondegenerate and irreducible surface $S$. An
$\alpha_2$--congruence $(\Sigma,p,q)$ consists of the tangent planes to a $1$--dimensional
family of cones which are a $1$--dimensional linear system of hyperplane sections of a
hypersurface in ${\bf P}^4$.
\end{corolario}

Finally, we give a proof of the footnote assertion due to Corrado Segre in [4] in
relation with our classification. We call Focal Variety to the projective realization
of the first order focal locus of the plane congruence $(\Sigma,p,q)$.

\begin{teorema}[Corrado Segre] Let $(\Sigma,p,q)$ be a
nondegenerate plane congruence in ${\bf P}^4$ such
that every first order focal point is a second order focal point.
\begin{enumerate}
\item If the general focal conic is nondegenerate, the focal variety
is a projection of the Veronese's surface, and the
congruence consists of the family of planes containing the
$2$--dimensional family of conics of the surface.
\item If the general focal conic consists of two different lines, then one of the next
possibilities is satisfied:
\begin{enumerate}
\item The focal variety consists of two cones in ${\bf
P}^4$, eventually coincident, with the same vertex $O$,
and the congruence consists of the planes generated by two generators, one in each
cone. This is a degenerate case of $\beta_3$--congruence.
\item The focal variety consists of a nondevelopable ruled
surface $S$ with plane directrix $C$ and the plane $\pi$
containing $C$. The congruence consists of the planes containing one generator
$r\subset S$ and one of the lines in $\pi$ passing through $r\cap \pi$. This is a
degenerate case of $\alpha_3$--congruence.
\end{enumerate}
\item If the general focal conic is a double line, the focal variety consists of:
\begin{itemize}
\item A line; the congruence consists of the family of planes
containing the line. A $\delta$--congruence.
\item A nondevelopable ruled surface; the congruence consists of the family of tangent
planes to that surface. This is a degenerate case of $\gamma_1$--congruence.
\item A developable ruled surface, which may be either a tangent
developable surface in the general case of $\gamma_2$--congruences or a cone, providing a
 degenerate case of 
$\gamma_3$--congruence; the congruence consists of linear 
$1$--dimensional families of planes passing through each one of the 
generators of the surface and containing the tangent plane along  
the generator. 
\end{itemize}
\end{enumerate}
\end{teorema}

{\bf Proof:} To prove property $1$, observe that the projective realization of the
family
${\bar q}:F_1(U)\longrightarrow U$ satisfies $\dim p(F_1(U))<\dim F_1(U)$, so the focal
variety is a surface containing a $2$--dimensional family of conics.

\smallskip

To prove property $2$, observe that if $C(\sigma)$ consists of two different lines
$r(\sigma)$ and
$r'(\sigma)$, the congruence is in case $\alpha$ or $\beta$. 

\smallskip

Suppose $(\Sigma,p,q)$ is a
$\beta$--congruence. By property $4$ in Theorem $2.1$, $p(R(U))$ and $p(R'(U))$ are
developable ruled surfaces. If one of them, for example $p(R(U))$, were a tangent
developable, by the classification in [3], $2.2.2$, $p(R(U))$ is the tangent developable
to the curve $\Sigma'$ consisting of the singular points $P(\sigma)$ of the focal
conics. Any generator $s$ of
$p(R'(U))$ meets every generator of $p(R(U))$ in the point $P(\sigma)$ and then
$\Sigma'$ is a line; in fact, $\Sigma'$ is the generator $s$. This is false, so $p(R(U))$ and
$p(R'(U))$ are cones with the same vertex $\Sigma'$ and this proves $(a)$.

\smallskip

Suppose $(\Sigma,p,q)$ is an $\alpha$--congruence and $r(\sigma)$ corresponds to the
developable direction. By Theorem $2.1$, the projective realization of these lines,
$p(R(U))$, is a non developable ruled surface, and by Theorem $2.2$, $R'(U)$ consists of
the lines in the plane $\pi=\overline{p(R'(U))}$. Finally, since every generator of
$p(R(U))$ meets a
$1$--dimensional family of lines in $\pi$, $p(R(U))$ has a plane directrix $C\subset\pi$.
This proves $(b)$.

\smallskip

To prove  property $3$, observe that the focal conic is $C(\sigma)=(r(\sigma))^2$, so
the congruence is in cases $\delta$ or $\gamma$. The projective realization of the
lines
$r(\sigma)$ consists of either a line in case $\delta$ or a ruled surface in 
case $\gamma$. We finish by applying properties $1$ and $3$ of Theorem
$2.1$. \rule{2mm}{2mm}

\end{document}